\documentclass[12pt,oneside]{amsart}
\usepackage{amssymb}
\usepackage[left=1in,top=1in,right=1in,bottom=1in]{geometry}
\newcommand{\ip}[2]{\langle#1,#2 \rangle}
\newcommand{\Co}{\mathrm{Co}} \newcommand{\cdual}{{\ast}}
\newcommand{\supp}{\mathrm{supp}}
\renewcommand{\theenumi}{\alph{enumi}} 
\theoremstyle{definition}
\newtheorem{theorem}{Theorem}[section]
\newtheorem{lemma}[theorem]{Lemma}
\newtheorem{assumption}[theorem]{Assumption}
\newtheorem{corollary}[theorem]{Corollary}
\newtheorem{proposition}[theorem]{Proposition}
\newtheorem{definition}[theorem]{Definition}
\newtheorem{example}[theorem]{Example}
\newtheorem{remark}[theorem]{Remark}

\begin{document}

\title{Examples of Coorbit Spaces for Dual Pairs} 

\author{Jens Gerlach Christensen} 
\address{329 Lockett Hall, Department of mathematics, 
  Louisiana State University} 
\email{vepjan@math.lsu.edu}
\urladdr{http://www.math.lsu.edu/\textasciitilde vepjan}
\thanks{The first author gratefully acknowledges support from the
  Louisiana Board of Regents under grant LEQSF(2005-2007)-ENH-TR-21
  and NSF grant DMS-0801010}

\author{Gestur \'Olafsson}
\address{322 Lockett Hall, Department of mathematics, 
  Louisiana State University} 
\email{olafsson@math.lsu.edu}
\urladdr{http://www.math.lsu.edu/\textasciitilde olafsson}
\thanks{The research of the second author was supported by
  NSF grants DMS-0402068 and DMS-0801010}

\begin{abstract}
  In this paper we summarize and give examples of
  a generalization of the coorbit space theory
  initiated in the 1980's by H.G. Feichtinger and K.H. Gr\"ochenig.
  Coorbit theory has been a powerful tool in characterizing Banach spaces
  of distributions with the use of integrable representations of locally
  compact groups. Examples are a wavelet characterization of the Besov
  spaces and a characterization of some Bergman spaces by the discrete
  series representation of $\mathrm{SL}_2(\mathbb{R})$. We present
  examples of Banach spaces which could not be covered by the
  previous theory, and we also provide atomic decompositions for
  an example related to a non-integrable representation.
\end{abstract}

\maketitle
\begin{center} {\today}
\end{center}

\section{Introduction}
\label{sec:introduction}
\noindent
In the 1980's H. G. Feichtinger and K. H. Gr\"ochenig developed a
unified theory for construction of Banach spaces and atomic
decompositions using representation theory (see \cite{FG0,FG1,FG2}).  
The construction, called coorbit space theory,
relies on an irreducible unitary representation $(\pi,\mathcal{H})$ of
a locally compact Hausdorff group $G$ with left invariant Haar measure $dx$. A
further requirement is that $\pi$ is integrable, in other words it is
required that there exists a non-zero $u\in \mathcal{H}$ such that
\begin{equation*}
  \int_G |(\pi(x)u,u)_\mathcal{H}|\, dx < \infty
\end{equation*}
Other constructions carry similar features as the theory of
Feichtinger and Gr\"ochenig, though they are not related to an
integrable nor irreducible representation.  For example the Banach
spaces of band-limited functions on a locally compact abelian group $G$ 
in \cite{FeiPan}, are related to the regular representation $\ell$ on
$L^2(G)$. This representation is neither irreducible nor integrable 
(for non-trivial $G$) and thus the spaces of band-limited function 
cannot be described 
as coorbit spaces. In \cite[Section 7.3]{FG0} it is also suggested
that coorbit spaces related to
non-integrable representations seem possible. 
The authors have developed a more general theory in \cite{CO}, which
does neither require integrability nor irreducibility. 
In the present article we give a short description of this theory,
but the main aim is to give several examples explaining
the general construction. Most of these examples are
not covered by the theory of Feichtinger and Gr\"ochenig. 

In section~\ref{sec:coorbitgen} we present our generalization of
the coorbit space theory of Feichtinger and Gr\"ochenig.
Next we show why smooth representations will make the construction
easier and we conclude the article with examples, some of which illustrate
the use of the generalized coorbit space theory. The examples include
coorbit spaces for Sobolev spaces on the $(ax+b)$-group, Bergman
spaces and discrete series representation of
$\mathrm{SL}_2(\mathbb{R})$, and bandlimited functions on a
Gelfand pair $(G,K)$ where $K$ is compact.

\section{Generalized Coorbit spaces}
\label{sec:coorbitgen}
\noindent
In this section we will suggest a minimal set of requirements for the
construction of general coorbit spaces.  This minimal set of
requirements will allow us to lift the restriction of only working
with integrable representations, and instead work with the bigger
class of square integrable representations and also non-irreducible
representations. 
Most proofs have been omitted, and we refer the reader to \cite{CO}
for the details.

Let $S$ be a Fr\'echet space and let $S^\cdual$ be the space of
continuous conjugate linear functionals on $S$ 
equipped with the weak topology. 
We assume that $S$ is
continuously imbedded and 
weakly dense in $S^\cdual$. The conjugate
dual pairing of elements $v\in S$ and $v'\in S^\cdual$ will be denoted
by $\ip{v'}{v}$.

Let $G$ be a locally compact group with a fixed left Haar measure
$dx$, and assume that $(\pi,S)$ is a representation of $G$. Also
assume that the representation is continous,i.e.  $g\mapsto \pi(g)v$
is continuous for all $v\in S$.  As usual define the contragradient
representation $(\pi^\cdual,S^\cdual)$ by
\begin{equation*}
  \ip{\pi^\cdual(x)v'}{v}=\ip{v'}{\pi(x^{-1})v}.
\end{equation*}
Then $\pi^*$ is a continuous representation of $G$ on $S^\cdual$. 
For a fixed vector 
$u\in S$ define the linear map $V_u:S^*\to C(G)$ by
\begin{equation*}
  V_u(v')(x) = \ip{v'}{\pi(x)u}.
\end{equation*}
The map $V_u$ is called \emph{the voice transform} or 
\emph{the wavelet transform}.
\begin{assumption}\label{assumption1}
  Let $Y$ be a left invariant Banach Space of Functions on $G$, and
  assume that there is a non-zero cyclic vector $u\in S$ 
  satisfying the following properties
  \begin{enumerate}
    \renewcommand{\theenumi}{R\arabic{enumi}}
  \item the reproducing formula $V_{u}(v)*V_{u}(u)=V_{u}(v)$ is true
    for all $v\in S$ \label{r1}
  \item the space $Y$ is stable under convolution with $V_u(u)$
    and $f \mapsto f*V_u(u)$ is continuous \label{r2}
  \item if $f=f*V_u(u)\in Y$ then the mapping
    $S \ni v \mapsto \int f(x) \ip{\pi^\cdual(x)u}{v} \,dx
    \in\mathbb{C}$ 
    is in $S^\cdual$\label{r3}
  \item the mapping 
    $S^\cdual\ni v' \mapsto \int \ip{v'}{\pi(x)u} 
    \ip{\pi^\cdual(x)u}{u}  \,dx \in\mathbb{C}$ 
    is weakly continuous \label{r4}
  \end{enumerate}
\end{assumption}
A vector $u$ satisfying Assumption~\ref{assumption1} 
is called an \emph{analyzing
vector}.  Note that (\ref{r4}) implies that there is an element
$v\in S$ such that
\begin{equation*}
  \ip{v'}{v} = \int \ip{v'}{\pi(x)u}\ip{\pi^\cdual(x)u}{u} \,dx
\end{equation*}
for all $v'\in S^\cdual$. If we use the notation $f^\vee(x) =
f(x^{-1})$ then $v = \pi(V_u(u)^\vee)u$.

\begin{theorem} \label{mainthm} Assume that $Y$ and $u$ satisfy
  Assumption~\ref{assumption1} and define the coorbit space
  \begin{equation}
    \Co_S^u Y = \{ v' \in S^\cdual | V_u(v') \in Y \}
  \end{equation}
  equipped with the norm $\| v' \| = \| V_u(v')\|_Y$.  Then the
  following properties hold
  \begin{enumerate}
  \item $V_u(v)*V_u(u) = V_u(v)$ for $v\in \Co^u_S Y$. \label{prop1}
  \item The space $\Co_S^u Y$ is a $\pi^\cdual$-invariant Banach
    space. \label{prop2}
  \item $V_u:\Co_S^u Y\to Y$ intertwines $\pi^\cdual$ and left
    translation \label{prop3}
  \item The convolution operator $f\mapsto f*V_u$ is a bounded
    projection from $Y$ to the closed subspace $V_u(\Co_S^u
    Y)=Y*V_u(u)$. \label{prop4}
  \item $\Co_S^u Y = \{ \pi^\cdual(f*V_u(u))u | f\in
    Y\}$. \label{prop5}
  \item $V_u:\Co_S^u Y \to Y*V_u(u)$ is an isometric
    isomorphism\label{prop6}
  \end{enumerate}
\end{theorem}

The following theorem tells us which analyzing vectors will give the
same coorbit space.
\begin{theorem}\label{thm:indep-of-u}
  If $u_1$ and $u_2$ both satisfy Assumption~\ref{assumption1} and for
  $i,j\in\{1,2\}$ the following properties can be verified
  \begin{itemize}
  \item there are non-zero constants $c_{i,j}$ such that
    $V_{u_i}(v)*V_{u_j}(u_i) = c_{i,j}V_{u_j}(v)$     for all $v\in S$
  \item $Y\ni f\mapsto f*V_{u_i}(u_j)\in Y$ is continuous
  \item $S^\cdual\ni v' \mapsto \int \ip{v'}{\pi(x)u_i}
    \ip{\pi^\cdual(x)u_i}{u_j} \,dx \in\mathbb{C}$ is weakly
    continuous
  \end{itemize}
  then $\Co_S^{u_1} Y = \Co_S^{u_2}Y$ with equivalent norms.
\end{theorem}

In the following we will describe how the choice of the Fr\'echet
space $S$ affects the coorbit space. We will show that there is great
freedom when choosing $S$.

\begin{theorem}
  Let $S$ and $T$ be Fr\'echet spaces which are weakly dense in their
  conjugate duals $S^\cdual$ and $T^\cdual$ respectively.  Let $\pi$
  and $\widetilde{\pi}$ denote representations of $G$ on $S$ and $T$
  respectively. Assume there is a vector $u\in S$ and
  $\widetilde{u}\in T$ such that the requirements in 
  Assumption~\ref{assumption1}
  are satisfied by both $(u,S)$ and
  $(\widetilde{u},T)$. Also assume that the conjugate dual
  pairings of $S^\cdual\times S$ and $T^\cdual\times T$ satisfy
  $\ip{u}{\pi(x)u}_S
  =\ip{\widetilde{u}}{\widetilde{\pi}(x)\widetilde{u}}_T$ for all
  $x\in G $.  Then $\Co_S^u Y$ and $\Co_T^{\widetilde{u}} Y$ are
  isometrically isomorphic. 
  The isomorphism is given by $V_{\widetilde{u}}V_u^{-1}$.

\end{theorem}

Let $\pi$ be a unitary irreducible representation of $G$ on
$\mathcal{H}$.  Assume that the Fr\'echet spaces $S$ and $T$ are
$\pi$-invariant and that $(S,\mathcal{H},S^\cdual)$ and
$(T,\mathcal{H},T^\cdual)$ are Gelfand triples with the common Hilbert
space $\mathcal{H}$.  Then $S\cap T$ is $\pi$-invariant and if we can
pick a non-zero vector $u\in S\cap T$, such that $u$ is analyzing for
both $S$ and $T$, then
\begin{equation*}
  \ip{u}{\pi(x)u}_S = (u,\pi(x)u)_{\mathcal{H}} = \ip{u}{\pi(x)u}_T
\end{equation*}
and we are in the situation of the previous theorem.  We summarize the
statement as
\begin{corollary}
  Assume that $(S,\mathcal{H},S^\cdual)$ and
  $(T,\mathcal{H},T^\cdual)$ are Gelfand tripples and assume there is
  an analyzing vector $u\in S\cap T$ such that both $(u,S)$ and
  $(u,T)$ satisfy Assumption~\ref{assumption1} for some Banach space
  $Y$, then $\Co_S^u Y$ and $\Co_T^{u} Y$ are isometrically isomorphic.
\end{corollary}

If the Fr\'echet space $S$ is a dense subspace of the Fr\'echet space $T$,
and $S$ is continuously included in $T$, then we can regard the space
$T^\cdual$ as a subspace of $S^\cdual$.  With this identification the
two coorbit spaces will be equal. We state the following
\begin{theorem}
  Let $(\pi,\mathcal{H})$ be a unitary irreducible representation of
  $G$, and let $(S,\mathcal{H},S^\cdual)$ and
  $(T,\mathcal{H},T^\cdual)$ be Gelfand triples for which $(\pi,S)$
  and $(\pi,T)$ are representations of $G$.  Assume that $i:S\to T$ is
  a continuous linear inclusion and that there is $u\in S$ such that
  both $(u,S)$ and $(i(u),T)$ satisfy
  Assumption~\ref{assumption1}.  Then the map $i^\cdual$ restricted to
  $\Co_T^{i(u)} Y$ is an isometric isomorphism between $\Co_T^{i(u)}
  Y$ and $\Co_S^u Y$.
\end{theorem}

\begin{remark}
  If $(\pi,S)$ is a representation of $G$ and $u$ is a cyclic vector
  for which it is true that
  $\ip{\pi^*(x)u}{u} = \overline{\ip{u}{\pi(x)u}}$ for all
  $x\in G$ and both (\ref{r1}) and (\ref{r4}) are satisfied, then
  $\ip{v}{w}$ is an inner product on $S$.  The completion
  $\mathcal{H}$ of $S$ with respect to the norm $\| v\|_\mathcal{H} =
  \sqrt{\ip{v}{v}}$ is a Hilbert space. The representation $\pi$ will then
  extend to a unitary representation $\widetilde{\pi}$ on
  $\mathcal{H}$, but we will not be able to conclude that
  $\widetilde{\pi}$ is irreducible.  Therefore the construction of
  coorbit spaces also works for non-irreducible representations, as
  long as we choose a cyclic vector in the Fr\'echet space $S$.

  Note also that a reproducing formula has been constructed from
  a non-unitary representation in \cite{Zimmermann}, thus
  allowing for construction of coorbit spaces in this new setting.
\end{remark}

\section{Smooth square integrable representations}
\label{sec:smooth-repr}
\noindent
We will now point out why smooth representations are easy to work
with. Let $(\pi,\mathcal{H})$ be a square integrable representation of
a Lie group $G$. This means that $\pi$ is unitary and irreducible
and that there is a non-zero $u$ such that
\begin{equation*}
  \int_G |(u,\pi(x)u)_\mathcal{H}|^2 \,dx <\infty
\end{equation*}
For square integrable representations the following theorem
found in \cite{DM} is crucial
\begin{theorem}[Duflo-Moore]
  \label{duflomoore}
  If $(\pi,\mathcal{H})$ is square integrable, then there is a positive
  densely defined operator $C$ with domain $D(C)$ such that
  $\ip{u}{\pi(\cdot)u}$ is in $L^2(G)$ if and only if $u\in D(C)$.
  Furthermore, if $u_1,u_2 \in D(C)$ then
  \begin{equation*}
    \int_G (v_1,\pi(x)u_1)_\mathcal{H}(\pi(x)u_2,v_2)_\mathcal{H} \, dx
    = (Cu_2,Cu_1)_\mathcal{H} (v_1,v_2)_\mathcal{H}
  \end{equation*}
\end{theorem}
Picking $u$ such that $\| Cu\| = 1$ we are able to obtain
the reproducing formula
\begin{equation}
  \label{eq:repsqint}
  V_u(v)*V_u(u) = V_u(v)
\end{equation}
for $v\in \mathcal{H}$. 
Denote the space of smooth vectors in $\mathcal{H}$  
by $\mathcal{H}_\pi^\infty$ equipped with
the usual Fr\'echet space topology (see \cite[p. 253]{warner1}).
Denote by 
$\mathcal{H}_\pi^{-\infty}$ the conjugate dual of
$\mathcal{H}_\pi^\infty$.  A vector $\pi(f)v = \int f(x)\pi(x)v \,dx$
for $f\in C_c^\infty(G)$ and $v\in\mathcal{H}$ is called a G{\aa}rding
vector.  In the following we let the Fr\'echet space $S$ be
$\mathcal{H}_\pi^\infty$. The reproducing formula
(\ref{eq:repsqint}) is then true for $v\in \mathcal{H}_\pi^\infty$,
and the following result is used to extend it
to $v\in \mathcal{H}_\pi^{-\infty}$.

\begin{proposition}\label{proposition:1}
  If $u\in \mathcal{H}_\pi^\infty$ is in the domain of the operator
  $C$ from Theorem \ref{duflomoore}, then the map
  \begin{equation*}
    \mathcal{H}_\pi^{-\infty} \ni v' \mapsto 
    \int \ip{v'}{\pi(x)u} \ip{\pi(x)u}{u}  \,dx \in \mathbb{C}
  \end{equation*}
  is continuous in the weak topology, and thus (\ref{r4}) is
  satisfied.
\end{proposition}

The next two results state, that if there is a non-zero $u\in
\mathcal{H}$ satisfying (\ref{r1}) and (\ref{r2}), then there is a
smooth non-zero vector which satisfies the same conditions.

\begin{proposition}
  Let $f\in C_c^\infty(G)$ then
  $\pi(f)D(C) \subseteq D(C)$.  
  If $u\in \mathcal{H}$ satisfies (\ref{r1}), then so
  will a constant multiple of $\pi(f)u$.
\end{proposition}

\begin{proof}
  Assume that $u\in D(C)$ and $f\in C_c^\infty (G)$.  Then
  \begin{align*}
    V_{\pi(f)u}(\pi(f)u)(z) &= \int\!\! \int f(x)
    f(y) (\pi(y)u,\pi(zx)u)_\mathcal{H} \,dx \,dy  \\
    &= f*V_u(u)*({f}^{\vee})(z)
  \end{align*}
  where ${f}^{\vee}(x) = f(x^{-1})$.
  Since $f,{f}^{\vee}\in C_c(G)\subseteq L^1(G)$ and $V_u(u)\in
  L^2(G)$, the function $f*V_u(u)*{f}^{\vee}$ is in $L^2(G)$.
  This shows that $\pi(f)u$ is in the domain of $C$.
\end{proof}

A Banach Function space $Y$ is called \emph{solid} if
$|f|\leq |g|$ a.e. and $g\in Y$ imply that $f\in Y$.

\begin{proposition}
  Assume that $Y$ is a solid Banach Function space, and that
  $Y\ni F\mapsto F*g \in Y$ is continuous for each $g\in
  C^\infty_c(G)$. Further assume that there is a $u\in \mathcal{H}$
  such that $Y\ni F\mapsto F*|V_u(u)|\in Y$ is continuous.  Then any
  G{\aa}rding vector of the form $\pi(f)u$ satisfies (\ref{r2}).
\end{proposition}

\begin{remark}
  Similar results can be shown for the space $\mathcal{H}_\pi^\omega$
  of analytic vectors, and its dual $\mathcal{H}_\pi^{-\omega}$, by
  approximating the G{\aa}rding vector $\pi(f)u$ by an analytic vector
  cf. \cite[p. 278 ff.]{warner1}.
\end{remark}

With these results in place the Gelfand triples
$(\mathcal{H}_\pi^\infty,\mathcal{H},\mathcal{H}_\pi^{-\infty})$ and
$(\mathcal{H}_\pi^\omega,\mathcal{H},\mathcal{H}_\pi^{-\omega})$ will
be natural choices for the construction of some coorbit spaces.

\section{Examples}
\label{sec:examples}
\noindent
In this section we will give some examples of coorbit spaces and also
show that the generalization presented in this paper is in fact a
generalization of the theory by Feichtinger and Gr\"ochenig.

\subsection{Coorbit spaces by Feichtinger and Gr\"ochenig}
\label{sec:coorb-spac-feicht}
\noindent
Here we will show that the coorbit spaces of Feichtinger and
Gr\"ochenig are special cases of the general coorbit spaces.
In the theory of Feichtinger and Gr\"ochenig (see \cite{FG0,FG1,FG2}
the Banach Function Space $Y$ is assumed to be left- and right-invariant
and also a right $L_m^1(G)$-module for some submultiplicative weight
$m$ satisfying $m(x)\geq 1$.  It is assumed that
\begin{equation*}
  \| f*g\|_Y \leq \|f\|_Y\|g\|_{L^1_m}.   
\end{equation*}
Let $(\pi,\mathcal{H})$ be a square-integrable irreducible unitary
representation and let
$$V_u(v)(x) = (v,\pi(x)u)_\mathcal{H}$$ 
and assume that
$$A_m = \{u\in \mathcal{H} | V_u(u)\in L_m^1(G) \}$$ 
is not trivial.  Then we can pick a non-zero $u\in A_m$ and define
\begin{equation*}
  \mathcal{H}_m^1=\{u\in \mathcal{H} | V_u(v)\in L_m^1(G) \}  
\end{equation*}
equipped with the norm $\| v\|_{\mathcal{H}_m^1}= \| V_u(v)\|_{L^1_m}$.
In \cite[Lemma 4.2]{FG0} Feichtinger and Gr\"ochenig show that
$\mathcal{H}_m^1$ is a $\pi$-invariant Banach space independent of the
chosen $u\in A_m$ (where different $u$ give equivalent norms). They also
show that $(\pi,\mathcal{H}_w^1)$ is irreducible and therefore every vector
$u\in A_w$ is cyclic.  These assumptions ensure that (\ref{r1}) and
(\ref{r2}) are automatically satisfied as has been proven in
\cite[Theorem 4.1]{FG1}.  By the same theorem the conjugate dual
$(\mathcal{H}_m^1)^\cdual$ of $\mathcal{H}_m^1$ consists of linear
functionals $v'$ for which $V_u(v')\in L_{1/m}^\infty(G)$, and thus
also (\ref{r4}) is satisfied.  It requires some work
(see \cite[Proposition 4.3]{FG1}) to show that if $f\in Y$ satisfies
$f=f*V_u(u)$ then $f\in L_{1/m}^\infty (G)$ and thus $\pi(f)u$ defined
by
\begin{equation*}
  \ip{\pi(f)u}{v} = \int f(x)\ip{\pi(x)u}{v} \,dx
\end{equation*}
is an element in $(\mathcal{H}_w^1)^\cdual$.  Therefore also (\ref{r3}) is
satisfied and the coorbit space $\Co Y$ of Feichtinger and Gr\"ochenig
is the space
\begin{equation*}
  \Co Y = \Co_{\mathcal{H}_m^1}^u Y = \{v'\in (\mathcal{H}_w^1)^\cdual | V_u(v')\in Y \}.
\end{equation*}
with norm
\begin{equation*}
  \| v'\|_{\Co Y} = \| V_u(v') \|_Y
\end{equation*}
The independence of $u$ follows from the construction of
$\mathcal{H}_m^1$ and Theorem~\ref{thm:indep-of-u}.

\subsection{Homogeneous Besov Spaces}
\label{sec:homog-besov-spac}
\noindent
We will now go through the calculations showing that
some homogeneous Besov spaces can be described as coorbits.
This gives a wavelet decomposition of these spaces. The example
has already been covered in \cite{FG0,FG1,Rauhut2}. 

The non-connected $(ax+b)$-group $G$ is the semidirect product
$\mathbb{R}^*\ltimes \mathbb{R}$, i.e.  the set
$\mathbb{R}^*\times \mathbb{R}$ with group composition
$(a,b)(a',b') = (aa',ab'+b)$.  The left Haar measure on the
$ax+b$-group is given by
\begin{equation*}
  f\mapsto \int_{\mathbb{R}}\int_{\mathbb{R}^*} f(a,b) \,\frac{da\,db}{a^2}
\end{equation*}
for $f\in C_c(G)$, where $da$ and $db$ are Lebesgue measures
on $\mathbb{R}$.  The representation
\begin{equation}
  \label{eq:repaxb}
  \pi(a,b)v(t) = \frac{1}{\sqrt{|a|}} v\Big(\frac{t-b}{a}\Big)
\end{equation}
is an irreducible unitary representation of $G$ on
$\mathcal{H}=L^2(\mathbb{R})$.  In the Fourier picture $\pi$ has the
form
\begin{equation}
  \label{eq:repaxb1}
  \widehat{\pi}(a,b)\widehat{v}(t) = e^{ibt} \sqrt{|a|} \widehat{v}(at)
\end{equation}

Let $\mathcal{S}(\mathbb{R})$ be the space of rapidly decreasing smooth
functions with the usual topology, and let $\mathcal{S}_0$ be the space
\begin{equation*}
  \mathcal{S}_0= \Big\{v \in \mathcal{S}(\mathbb{R}) 
  \Big| \int x^nv(x) \,dx= iD^n\widehat{v}(0)= 0  \Big\}  
\end{equation*}
with the subspace topology. $\mathcal{S}_0$ is a $\pi$-invariant
Fr\'echet space. The conjugate dual space $\mathcal{S}_0^\cdual$ is
isomorphic to
\begin{equation*}
  \mathcal{S}_0^\cdual \cong \mathcal{S}'(\mathbb{R})/\{\Lambda_P|\text{$P$ is a polynomial} \}  
\end{equation*}
where $\Lambda_P$ is the functional given by $\Lambda_P (v) = \int
P(x)v(x)\,dx$ for $v\in \mathcal{S}(\mathbb{R})$.  Denote the
conjugate dual pairing by $\ip{v'}{v}$ (which extends the inner
product on $L^2(\mathbb{R})$) and the voice transform
\begin{equation*}
  V_v(v')(x) = \ip{v'}{\pi(x)v}
\end{equation*}
for $v\in \mathcal{S}_0$ and $v'\in \mathcal{S}_0'$.

It is known (see \cite[Theorem 14.0.2, p. 66]{holschneider}) that if
\begin{equation}
  \label{eq:admissible}
  0 < \int_{-\infty}^\infty \frac{|\widehat{u}(x)|^2}{|x|} \,dx < \infty
\end{equation}
then for any $v\in \mathcal{S}_0$
the inversion formula for the wavelet transform gives
\begin{equation*}
  v = C_{u} \int_{\mathbb{R}} 
  \int_\mathbb{R} V_u(v)(a,b) \pi(a,b)u \,\frac{da\,db}{a^2}
\end{equation*}
for some constant $C_u$ depending only on $u$. Since all non-zero
functions in $\mathcal{S}_0$ satisfy (\ref{eq:admissible}) they are
all cyclic in $S_0$.
We see that $L^2(\mathbb{R})$ is
weakly dense in $\mathcal{S}'_0$,
and thus $\mathcal{S}_0$, which is
dense in $L^2(\mathbb{R})$ by irreducibility, is also weakly dense in
$\mathcal{S}_0'$.

To check the remaining properties needed for construction of coorbit
spaces, let us note the following result which is equivalent to
\cite[Theorem 19.0.1]{holschneider}
\begin{lemma} \label{lem:holschneider} Let
  $$
  \mathcal{S}(\mathbb{R}_0\times \mathbb{R}) = \{v\in
  \mathcal{S}(\mathbb{R}^2) | \sup|(a^2+a^{-2})^n(1+b^2)^m D^\alpha
  \widehat{v}(a,b)| < \infty \}
  $$
  then
  $$
  V_u: \mathcal{S}_0 \to \mathcal{S}(\mathbb{R}_0\times \mathbb{R})
  $$
  is continuous for any $u\in \mathcal{S}_0$.
\end{lemma}

(\ref{r1}) By Lemma~\ref{lem:holschneider}
we see that $V_u(v)$ is in $L^1(G)$ and
thus also in $L^2(G)$ when $u,v\in \mathcal{S}_0(\mathbb{R})$.
Therefore we can pick a $u\in \mathcal{S}_0$ such that $V_u(v)*V_u(u)
= V_u(v)$ for all $v\in \mathcal{S}_0$.

(\ref{r2}) Let $Y$ be the mixed norm Banach space 
\begin{equation*}
  L_s^{p,q}(G) 
  = \Big\{f\,\Big|\, \Big( \int_{\mathbb{R}^*} 
  \Big( \int_{\mathbb{R}} |f(a,b)|^p a^{-s}
  \,db \Big)^{/p} \frac{da}{a^2} \Big)^{1/q} < \infty
  \Big\}
\end{equation*}
for $p,q \geq 1$ and $s\in\mathbb{R}$.
Then by
Minkowski's inequality
\begin{equation*}
  \| f*V_u(u)\|_{L_s^{p,q}} \leq \| f\|_{L_s^{p,q}} \| V_u(u)\|_{L_s^1}
\end{equation*}
which shows that $f\mapsto f*V_u(u)$ is continuous.

(\ref{r3}) From Lemma~\ref{lem:holschneider} it is obtained that for
$1\leq p',q'\leq \infty$ the voice transform $V_u(u)$ is in
$L_{-s}^{p',q'}(G)$.
Then H\"olders inequality tells us that $f*V_u(u) \in L_{-s}^\infty(G)$ if
we pick $p'$ and $q'$ such that $1/p+1/p' =1$ and $1/q +1/q'=1$.  So
if $f=f*V_u(u)$ it will follow that $f\in L_{-s}^\infty(G)$ and then
\begin{equation*}
  \mathcal{S}_0 \ni v\mapsto 
  \int_\mathbb{R}\int_\mathbb{R} f(a,b) \ip{v}{\pi(a,b)u} \,\frac{da\,db}{a^2}
\end{equation*}
is continuous, since the mapping $S_0\ni v\mapsto V_u(v) \in L_s^1(G)$
is continuous.

(\ref{r4}) That
\begin{equation*}
  \mathcal{S}_0' \ni v' \mapsto
  \int_\mathbb{R}\int_\mathbb{R} 
  \ip{v'}{\pi(a,b)u}\ip{\pi(a,b)u}{u} \,\frac{da\,db}{a^2}
\end{equation*}
is weakly continuous is exactly the statement of Theorem 24.1.1 in
\cite{holschneider}.

This verifies that the space
\begin{equation*}
  \Co_{\mathcal{S}_0}^u L_s^{p,q}(G) 
  = \{v'\in \mathcal{S}_0' | V_u(v) \in L_s^{p,q}(G) \}
\end{equation*}
is Banach a space. Since the representation $\pi$ is integrable this
space coincides with the coorbit space $\Co_{FG} L_s^{p,q}(G)$ of
Feichtinger and Gr\"ochenig. In \cite{FG0,Grochenig1} it is proven that
the homogeneous Besov spaces (see \cite{Triebel1988b})
\begin{equation*}
  \dot{B}_{p,q}^s = \Big \{f\in S_0'\Big| \sum_{j\in\mathbb{Z}} 2^{sjq} \|
  \mathcal{F}^{-1}(\phi_j\mathcal{F}f) \|^q_{L^p} < \infty   \Big\}
\end{equation*}
can be characterized as $\Co_{\mathcal{S}_0}^u L_{s+1/2-1/q}^{p,q}(G)
$. Here $\phi\in \mathcal{S}_0$ and $\phi_j(x)=\phi(2^{-j}x)$
are chosen such that 
$\supp(\widehat{\phi}) \subseteq \{ x\, |\, 1/2 \leq |x| \leq 2  \}$
and $\sum_{j\in\mathbb{Z}} \phi_j =1$.

\subsection{Coorbit spaces for Sobolev spaces on the
  \boldmath${(ax+b)}$-group}
\label{sec:coorb-spac-sobol}
\noindent 
We now present a type of coorbit spaces which cannot
be described by previous coorbit theories. In particular
we find coorbit spaces related to the non-solid 
Sobolev spaces on the $(ax+b)$-group.

Again let $G$ be the $(ax+b)$-group with representation 
$(\pi,\mathcal{H})$ as before.
The $(ax+b)$-group can also be regarded as a matrix group
\begin{equation*}
  G \simeq \Big\{ 
  \begin{pmatrix}
    a & b \\ 0 & 1
  \end{pmatrix} \Big| a\neq 0, b\in\mathbb{R}
  \Big\}
\end{equation*}
The matrices 
\begin{equation*}
  X_1 = 
  \begin{pmatrix}
    1 & 0 \\ 0 & 0
  \end{pmatrix}
\qquad\text{and}\qquad
  X_2 = 
  \begin{pmatrix}
    0 & 1 \\ 0 & 0
  \end{pmatrix}
\end{equation*}
form a basis for the Lie algebra $\mathfrak{g}$.
The representation $\pi$ induces the following
differential operators
\begin{equation*}
  \pi(X_1)f(x) = -\frac{1}{2}f(x) - xf'(x) \qquad\text{and}\qquad
  \pi(X_2)f(x) = -f'(x).
\end{equation*}
For $f\in C^\infty(G)$ and $X\in\mathfrak{g}$ let
\begin{equation*}
Xf(x) = \frac{d}{dt}\Big|_{t=0} f(\exp(-tX)x).  
\end{equation*}
For $p\geq 1$ define the Sobolev space
(see \cite[Definition 3]{Triebel1})
\begin{equation*}
  W^m_p(G)
  = \Big\{ 
  f \in L^p(G) \,|\, 
  \| f\|_{W^m_p} = \sum_{k=1}^m \sum_{n_k\in \{ 1,2\} } \| X_{n_1}\cdots X_{n_k} f \|_{L^p} 
  < \infty  \Big\}
\end{equation*}
As in the previous section we use the Gelfand triple
$(\mathcal{S}_0,\mathcal{H},\mathcal{S}_0^*)$ for which the conditions
(\ref{r1}) and (\ref{r4}) are satisfied.  Since $W^m_p(G) \subseteq
L^p(G)$ and $V_u(v)\in (L^p(G))^*$ for $p\geq 1$ the condition
(\ref{r3}) is also satisfied.  It remains to show that
\begin{equation*}
  W^m_p(G) \ni f\mapsto f*V_u(u) \in W^m_p(G)  
\end{equation*}
is a continuous mapping for some $u\in \mathcal{S}_0$.  We note again
that by Lemma~\ref{lem:holschneider}
$V_u(v) \in L^1(G)$ for any $u,v\in \mathcal{S}_0$
and therefore $f*V_u(u)\in L^p(G)$ for $f\in W^m_p(G)$. 
Finally 
\begin{equation*}
  X_i (f*V_u(u))
  = f*V_{\pi(X_i)u}(u),
\end{equation*}
and since $\pi(X_i)\mathcal{S}_0\subseteq \mathcal{S}_0$
it follows that $X_i(f*V_u(u))\in L^p(G)$.  
Thus $f*V_u(u)\in W^m_p(G)$ if $f\in L^p(G)$, and
$W^m_p(G)\ni f\mapsto f*V_u(u) \in W^m_p(G)$ is continuous.
Therefore the space
\begin{equation*}
  \Co_{\mathcal{S}_0}^u W^m_p(G)
  = \{ v'\in \mathcal{S}_0 | V_u(v')\in W^m_p(G)  \}
\end{equation*}
is a Banach space.  Notice that $W^m_p(G)$ is not solid for $m\geq 1$
(the characteristic functions for compact intervals are not in
$W^m_p(G)$).

\subsection{Discrete series representation for
  \boldmath$\mathrm{SU}(1,1)$ and the \boldmath$(ax+b)$-group}
\label{sec:discr-seri-repr}
In this section we show that some Bergman spaces
can be described as coorbit spaces related to the
discrete series representations $\pi_n$ of 
$\mathrm{SL}_2(\mathbb{R})$ for integers $n\geq 2$. 
The construction will also work for representations
$\pi_s$ for non-integer parameters $s>1$ if we
replace $SL_2(\mathbb{R})$ by a covering group.
This has been done
before in \cite{FG0} in the case of the integrable
representations. In that paper it is also noted,
that it should be possible to define coorbits for
a non-integrable discrete series representation.
With the generalized coorbit space theory, this is 
indeed possible as we will show in the following.
Furthermore we also show that the discretization
procedures of \cite{FG1,Grochenig1} 
carry over to this non-integrable case.

The example can be generalized
to other bounded symmetric domains, but that would 
require too much new notation to carry out here.

\subsubsection{Continuous description}
\noindent
The connected $(ax+b)$-group $G$ can be realized as the group
\begin{equation*}
  G = \Big\{ 
  \begin{pmatrix}
    a & b \\ 0 & a^{-1}
  \end{pmatrix}
  \Big| a>0, b\in\mathbb{R} \Big\}
  \subseteq \mathrm{SL}_2(\mathbb{R})
\end{equation*}
which can be regarded as a subgroup of $\mathrm{SU}(1,1)$ using the
Cayley transform:
\begin{align*}
  \begin{pmatrix}
    \alpha & \beta \\ \bar\beta & \bar\alpha
  \end{pmatrix}
  &= \frac{1}{2}
  \begin{pmatrix}
    1 & -i \\ -i & 1
  \end{pmatrix}
  \begin{pmatrix}
    a & b \\ 0 & a^{-1}
  \end{pmatrix}
  \begin{pmatrix}
    1 & i \\ i & 1
  \end{pmatrix} \\
  &= \frac{1}{2}
  \begin{pmatrix}
    a+a^{-1} +ib & b+i(a-a^{-1})\\ b-i(a-a^{-1}) & a+a^{-1} - ib
  \end{pmatrix}
\end{align*}

Let $\mathbb{D}$ be the unit disc $\mathbb{D}=\{z\in\mathbb{C}||z|<1
\}$ and $n\geq 2$ an integer, then the pairing
\begin{equation*}
  (u,v)_n 
  = \int_\mathbb{D} u(z) \overline{v(z)} (1-|z|^2)^{n-2} \,dx\, dy
\end{equation*}
gives and inner product on the Hilbert space
\begin{equation*}
  \mathcal{H}_n = \mathcal{A}^{2,n}(\mathbb{D}) 
  = \{ v\in \mathcal{O}(\mathbb{D})| (v,v)_n < \infty \}.
\end{equation*}
The discrete series representation
\begin{equation*}
  \pi_n
  \begin{pmatrix}
    \alpha & \beta \\ \bar\beta & \bar\alpha
  \end{pmatrix} v(z) = (-\bar\beta z + \bar\alpha)^{-n}
  v\Big(\frac{\alpha z -\beta}{-\bar\beta z +\bar\alpha }\Big)
\end{equation*}
is an irreducible unitary representation of $G$ on
$\mathcal{H}_n=\mathcal{A}^{2,n}(\mathbb{D})$ (see
\cite{Harish-Chandra1955,Harish-Chandra1956a,knapp,olafssonoersted}). 

As distribution space $S^\cdual$ we will use the 
conjugate dual of the
space $S=\mathcal{H}_n^\infty$.
We will denote the conjugate dual
by $\mathcal{H}_n^{-\infty}$.

The function $\psi(z) = 1$ is in $\mathcal{H}^\infty_n$ for all
$n$. Since $(\pi_n,\mathcal{H}_n^\infty)$ is irreducible 
as a representation $SL_2(\mathbb{R})$ (see \cite[p. 254]{warner1}),
the function $\psi$ is 
$\mathrm{SL}_2(\mathbb{R})$-cyclic in $\mathcal{H}_n^\infty$.
The group $SU(1,1)$ has the Iwasawa decomposition
$NAK$ where 
\begin{align*}
  N &= \Big\{
  \begin{pmatrix}
    1+ir & r \\ r &1-ir
  \end{pmatrix}
  \Big| r\in\mathbb{R} \Big\}, \\
  A &= \Big\{
  \begin{pmatrix}
    \cosh s & i\sinh s \\ -i\sinh s &\cosh s
  \end{pmatrix}
  \Big| s\in\mathbb{R} \Big\} \\
  K &= \Big\{
  \begin{pmatrix}
    e^{it} & 0 \\ 0 & e^{-it}
  \end{pmatrix}
  \Big| t\in\mathbb{R} \Big\},
\end{align*}
The action of an element
in $K$ on $\psi$ is
\begin{equation*}
  \pi_n
  \begin{pmatrix}
    e^{it} & 0 \\ 0 & e^{-it}
  \end{pmatrix}
  \psi (z)
  = e^{int} \psi(z),
\end{equation*}
and hence $\psi$ is also $NA$-cyclic in $\mathcal{H}_n^\infty$.
The group $NA$ is exactly the subgroup of $\mathrm{SU}(1,1)$
corresponding to $G$, so we conclude that $\psi$ is $G$-cyclic
in $\mathcal{H}_n^\infty$.

By use of polar coordinates and Cauchy's integral theorem
the voice transform $V_\psi(\psi)$ is
\begin{equation*}
  V_{\psi}(\psi) 
  \begin{pmatrix}
    \alpha & \beta \\ \bar\beta & \bar\alpha
  \end{pmatrix}
  = \int_\mathbb{D} (-\bar\beta z + \bar\alpha)^{-n} (1-|z|^2)^{n-2}
  \,dx\, dy
  = \frac{1}{2(n-1) \bar\alpha^{n}}
\end{equation*}

As a function of $(a,b)$ we have
\begin{equation*}
  F_n(a,b) 
  = V_{\psi}(\psi)  
  \begin{pmatrix}
    \alpha & \beta \\ \bar\beta & \bar\alpha
  \end{pmatrix}
  = \frac{1}{2(n-1)} (a+a^{-1} -ib)^{-n}e^{-int}.
\end{equation*}
Furthermore 
\begin{equation*}
  \| F_n\|_{L^p(G)}^p 
  = \frac{1}{2^p(n-1)^p}\int_{-\infty}^\infty \int_0^\infty
  \frac{1}{((a+a ^{-1})^2 + b^2)^{np/2}} \,\frac{da \,db}{a^2} < \infty
\end{equation*}
if and only if $np>2$.
In particular we recover the well know fact 
that $F_n$ is square integrable when $n\geq 2$
and integrable for $n\geq 3$.

While $F_2$ is not in $L^1(G)$ it is very close in the following sense
\begin{lemma}\label{lem:int1}
  If $0 < \epsilon < 2$ then 
  \begin{equation*}
    \int_0^\infty \int_{-\infty}^\infty 
    |F_2(a,b)| a^\epsilon \,\frac{db\, da}{a^2} < \infty
  \end{equation*}
\end{lemma}

This can be used to show that the conditions for construction of
coorbit spaces are satisfied.

(\ref{r1}) is satisfied since $F_n\in L^2(G)$ for $n\geq 2$, showing
that the representations $\pi_n$ are square integrable. 
Let $u$ 
be a normalization of $\psi$ in order to make convolution with
$V_u(u)$ an idempotent.

(\ref{r2}) The mapping $f\mapsto f*F_n$ is bounded from $L^p(G)$ to
$L^p(G)$ for $p > 1$.  Using Minkowski's inequality it is easy to see
for $n\geq 3$, since then $F_n$ is in $L^1(G)$. We will show that it
is also true for $n=2$.  Assume that $f\in L^p(G)$ and pick $q$ such
that $1/p+1/q =1$, and look at
\begin{align*}
  \Big |\int\!\! \int &f(a,b) F_2((a,b)^{-1}(a_1,b_1))
  \,\frac{da\, db}{a^2}\Big|^p \\
  &\leq \Big(\int\!\! \int |f(a,b)| a^{-1/pq}
  |F_2((a,b)^{-1}(a_1,b_1))|^{1/p+1/q} a^{1/pq}
  \,\frac{da\, db}{a^2}\Big)^p \\
  &\leq \Big( \int\!\! \int |f(a,b)|^p a^{-1/q}
  |F_2((a,b)^{-1}(a_1,b_1))| \,\frac{da\, db}{a^2} \Big) \\
  &\qquad \times \Big( \int\!\! \int a^{1/p}
  |F_2((a,b)^{-1}(a_1,b_1))| \,\frac{da\, db}{a^2} \Big)^{p/q}
\end{align*}
By the unitarity of $\pi_n$ it is true that $|F_n((a,b))| =
|F_n((a,b)^{-1})|$ so the second integral becomes
\begin{equation}\label{eq:1}
  \int\!\! \int a^{1/p}
  |F_2((a_1,b_1)^{-1}(a,b))| \,\frac{da\, db}{a^2}
  = 
  \int\!\! \int (aa_1)^{1/p}
  |F_2((a,b))| \,\frac{da\, db}{a^2}
  \leq Ca_1^{1/p}
\end{equation}
where the have used the invariance of the measure and 
Lemma~\ref{lem:int1}. Thus we get
\begin{align*}
  \Big |\int\!\! \int &f(a,b) F_2((a,b)^{-1}(a_1,b_1))
  \,\frac{da \,db}{a^2}\Big|^p \\
  &\leq C \int\!\! \int |f(a,b)|^p a^{-1/q} |F_2((a,b)^{-1}(a_1,b_1))|
  \,\frac{da\, db}{a^2} a_1^{1/q}.
\end{align*}
Now we can calculate an estimate for the norm using Fubuni's theorem
and (\ref{eq:1}) again
\begin{align*}
  \| f*F_2 \|_p^p &\leq C \int\!\! \int\!\! \int\!\! \int |f(a,b)|^p
  a^{-1/q} |F_2((a,b)^{-1}(a_1,b_1))| \,\frac{da\, db}{a^2} a_1^{1/q}
  \frac{da_1 db_1}{a_1^2} \\
  &= C \int\!\! \int |f(a,b)|^p a^{-1/q} \int\!\! \int
  |F_2((a,b)^{-1}(a_1,b_1))| a_1^{1/q}
  \,\frac{da_1 \,db_1}{a_1^2} \,\frac{da \,db}{a^2}  \\
  &\leq
  C^2 \int\!\! \int |f(a,b)|^p a^{-1/q} a^{1/q} \,\frac{da\, db}{a^2} \\
  &= C^2 \| f\|_p^p
\end{align*}
This shows that $f\mapsto f*F_n$ is a continuous linear operator on
$L^p(G)$ for all $n\geq 2$.

(\ref{r3}) The smooth vectors for the representation $\pi_n$ have been
characterized in \cite{olafssonoersted} and more generally
in \cite{Chebli2004} to be
\begin{equation*}
  \mathcal{H}_n^\infty 
  = \Big\{ \sum_{k=0}^\infty a_kz^k \Big| \forall m:
  \sum_{k=0}^\infty |a_k|^2\frac{(n-1)!k!}{(n+k-1)!}(n(n-2)+2k^2)^{2m}
  < \infty \Big \}
\end{equation*}
with conjugate dual consisting of formal power series
\begin{equation*}
  \mathcal{H}_n^{-\infty} 
  = \Big\{ \sum_{k=0}^\infty b_kz^k \Big| \exists m:
  \sum_{k=1}^\infty |b_k|^2\frac{(n-1)!k!}{(n+k-1)!}(n(n-2)+2k^2)^{-2m}
  < \infty \Big\} 
\end{equation*}

Let $\phi$ be a smooth vector with expansion $\sum_{k=0}^\infty
a_kz^k$.  The Taylor series for $ \pi
\begin{pmatrix}
  \alpha & \beta \\
  \bar\beta & \bar\alpha
\end{pmatrix} 1(z) $ is
\begin{equation*}
  \pi
  \begin{pmatrix}
    \alpha & \beta \\
    \bar\beta & \bar\alpha
  \end{pmatrix} 1(z) = \alpha^{-n} \sum_{k=0}^\infty (\beta/\alpha)^k
  (-1)^k \frac{(n+k-1)!}{(n-1)!k!} z^k
\end{equation*}
Therefore, since $(z^k,z^k)_n= k!(n-2)!/(k+n-1)!$,
\begin{equation*}
  V_\psi(\phi)(a,b)
  = \frac{1}{(n-1)\bar\alpha^n}  
  \sum_{k=0}^\infty (-1)^k (\beta/\alpha)^k a_k
  = 2 F_n(a,b) \sum_{k=0}^\infty (-1)^k (\beta/\alpha)^k a_k
\end{equation*}
and it can estimated by
\begin{equation*}
  |V_1(\phi)(a,b)| 
  \leq 2 |F_n(a,b)| \sum_{k=0}^\infty |a_k|
\end{equation*}
since $|\alpha|>|\beta|$.  For any $m$ we obtain that
\begin{align*}
  \sum_{k=0}^\infty |a_k| &\leq |a_0| + \Big(\sum_{k=1}^\infty
  \frac{(n+k-1)!}{k!(n-1)!}
  (n(n-2)+2k^2)^{-2m}\Big)^{1/2}  \\
  &\qquad\times\Big(\sum_{k=1}^\infty |a_k|^2
  \frac{(n-1)!k!}{(n+k-1)!}  (n(n-2)+2k^2)^{2m} \Big)^{1/2}
\end{align*}
and by picking $m$ large enough the first sum converges.  Thus the
mapping $\phi\mapsto V_\psi(\phi)$ is continuous from
$\mathcal{H}_\pi^\infty$ to $L^q$ for all $q>1$.  So if $f\in L^p(G)$
for $p>1$ then the mapping $\phi\mapsto \int f(a,b) V_\psi (\phi)(a,b)
\,\frac{da\,db}{a^2}$ is continuous by H\"older's inequality. This
ensures the validity of (\ref{r3}).

(\ref{r4}) By Proposition~\ref{proposition:1} the condition (\ref{r4})
is automatically satisfied when $u$ is a constant multiple of $\psi\in
\mathcal{H}_n^\infty$.

This means that all the conditions in Assumption~\ref{assumption1}
are satisfied
for some $u\in \mathcal{H}_n^\infty$ (constant multiple of
$\psi(z)=1$) with the Banach space $Y=L^p(G)$. Therefore the space
\begin{equation*}
  \Co_{\mathcal{H}_n^\infty}^u L^p(G)
  = \{v'\in \mathcal{H}_n^{-\infty} | V_u(v')\in L^p(G) \}
\end{equation*}
is a Banach space. As pointed out in \cite[Section 7]{FG0} these
spaces are the Bergman spaces $\mathcal{A}^{p,pn/2}$ where
\begin{equation*}
  \mathcal{A}^{p,\alpha}(\mathbb{D}) = \Big\{f\in\mathcal{O}(\mathbb{D}  
  \Big| \int_\mathbb{D} |f(z)|^p (1-|z|^2)^{\alpha-2} dz < \infty      \Big\}.
\end{equation*}

\subsubsection{Discretization}

In the following we will only work with the non-integrable
voice transform $F_2=V_\psi(\psi)$ since the other cases are covered by
the theory of Feichtinger and Gr\"ochenig. The aim is to show
that discretizations of coorbit spaces are possible without integrable
representations.

The key to finding atomic decompositions will be the following
result
\begin{lemma}\label{lem:osc}
  For each $\epsilon > 0$ there is a neighbourhood 
  $U$ of the identity such
  that 
  \begin{equation*}
    \Big| \frac{F_2((a,b)(x,y))}{F_2(x,y)} -1 \Big| <\epsilon
  \end{equation*}
  for $(a,b)\in U$.
\end{lemma}

\begin{proof}
  \begin{align*}
    \Big| \frac{F_2((a,b)(x,y))}{F_2(x,y)} -1 \Big|^2
    &=  \Big| \frac{(ax +(ax)^{-1}-
      i(bx^{-1}+ay))^2 -
      (x+x^{-1}-iy)^2}{(x+x^{-1}-iy)^2}
    \Big|^2 \\
    &=  \Big| \frac{ax +(ax)^{-1}-
      i(bx^{-1}+ay) +
      (x+x^{-1}-iy)^2}{x+x^{-1}-iy}
    \Big|^2 \\
    &\qquad\times
    \Big| \frac{ax +(ax)^{-1}-
      i(bx^{-1}+ay) -
      (x+x^{-1}-iy)^2}{x+x^{-1}-iy}
    \Big|^2 \\
    &=  \Big|
    \frac{(1+a)x+(1+a^{-1}x^{-1} - i(bx^{-1}+(1+a)y}{x+x^{-1}-iy}
    \Big|^2     \\
    &\qquad\times
    \Big|
    \frac{(a-1)x+(a^{-1}-1)x^{-1} -i(bx^{-1}+(a-1)y)}{x+x^{-1}-iy}
    \Big|^2 \\
  \end{align*}
  Let 
  $1/\delta < a < \delta$ for $\delta >1$ and
  $-\gamma < b <\gamma$ for $\gamma >0$ and then 
  estimate each of the squares
  \begin{align*}
    \Big|
    &\frac{(1+a)x+(1+a^{-1}x^{-1} - i(bx^{-1}+(1+a)y}{x+x^{-1}-iy}
    \Big|^2 \\
    &=
    \frac{((1+a)x+(1+a^{-1}x^{-1})^2
      +
      ((bx^{-1}+(1+a)y)^2}{(x+x^{-1})^2+y^2} \\
    &\leq 
    \frac{((1+\delta)x+(1+\delta) x^{-1})^2}{(x+x^{-1})^2+y^2}
    + \frac{b^2x^{-2}}{(x+x^{-1})^2+y^2} 
    +\frac{(1+\delta)^2y^2}{(x+x^{-1})^2+y^2}
    +\frac{2(1+\delta)|b||y|x^{-1}}{(x+x^{-1})^2+y^2} \\
    &\leq 
    (1+\delta)^2 + b^2 +(1+\delta)^2
    + 2(1+\delta)|b| 
    \frac{|y|}{\sqrt{(x+x^{-1})^2+y^2}}
    \frac{x^{-1}}{\sqrt{(x+x^{-1})^2+y^2}} \\
    &\leq 
    (1+\delta)^2 + \gamma^2 +(1+\delta)^2
    + 2(1+\delta)\gamma 
  \end{align*}
  This we can get arbitrarily close to $8$ by picking $\delta$ and
  $\gamma$ close enough to $1$ and $0$ respectively.
  Similarily
  \begin{align*}
    \Big|
    &\frac{(a-1)x+(a^{-1}-1)x^{-1} - i(bx^{-1}+(a-1)y}{x+x^{-1}-iy}
    \Big|^2 \\
    &=
    \frac{((a-1)x+(a^{-1}-1)x^{-1})^2
      +
      ((bx^{-1}+(a-1)y)^2}{(x+x^{-1})^2+y^2} \\
    &\leq 
    \frac{((\delta-1)x+(\delta-1) x^{-1})^2}{(x+x^{-1})^2+y^2}
    + \frac{b^2x^{-2}}{(x+x^{-1})^2+y^2}
    +\frac{(a-1)^2y^2}{(x+x^{-1})^2+y^2}
    +\frac{2(\delta-1)|b||y|x^{-1}}{(x+x^{-1})^2+y^2} \\
    &\leq 
    (\delta-1)^2 + \gamma^2 +(\delta-1)^2
    + 2(\delta-1)\gamma 
  \end{align*}
  This last term can be made arbitrarily close to $0$,
  thus finishing the proof of the lemma.
\end{proof}
From this result follows easily 
\begin{corollary}\label{cor:1}
  There exist a neighbourhood $U$ of the identity
  and constants $C_1,C_2 > 0$ such that
  \begin{equation*}
    C_1 |F_2(x,y)| \leq |F_2((a,b)(x,y))| \leq C_2 |F_2(x,y)|
  \end{equation*}
  for all $(x,y)\in G$ with $(a,b)\in U$.
  These constants can be chosen arbitrarily close to $1$,
  by choosing $U$ small enough.
\end{corollary}

When discretizing the reproducing formula,
we need the following definition from \cite{FG1}:
\begin{definition}
  Let $V$ be a compact neighbourhood of the identity.
  The set of points $\{ x_i\}$ are said to be $V$-separated if
  the $x_iV$ are pairwise disjoint.
  
  Let $U$ be a compact neighbourhood of the identity.
  The set of points $\{ x_i\}$ are said to be $U$-dense 
  if $G=\cup_i x_iU$.
\end{definition}

\begin{proposition}
  Let $V\subseteq U$ be compact neighbourhoods of the identitiy.
  Assume that the points $\{ x_i\}$ are $V$-separated and $U$-dense
  and that $U$ satisfies Corollary~\ref{cor:1}. Let 
  $\{ \psi_i\}$ be a partition of unity 
  for which $\supp(\psi_i)\subseteq x_iU$.
  Then the following is true
  \begin{enumerate}
  \item The mapping 
    $\ell^p \ni (\lambda_i) \mapsto \sum_i \lambda_i \ell_{x_i} F_2
    \in L^p(G)*F_2$ is continuous \label{item:1}
  \item The mapping 
    $L^p(G)*F_2 \ni f \mapsto (f(x_i))_{i\in I} \in \ell^p(I)$
    is continuous\label{item:2}
  \item The mapping 
    $L^p(G)*F_2 \ni f\mapsto (\int_G f(x)\psi_i(x) dx)_{i\in I} \in
    \ell^p(I)$
    is continuous \label{item:3}
\end{enumerate}
\end{proposition}

\begin{proof}

  First note that the norms on $\ell^p$ and $L^p(G)$ are related
  in the following sense
  \begin{equation*}
    \| (\lambda_i) \|_{\ell^p} 
    = \frac{1}{|V|} \Big\| \sum_i \lambda_i 1_{x_iV} \Big\|_{L^p}
  \end{equation*}
  Also the convolution in $L^p(G)$ with $|F_2|$ is continuous and we will
  denote the norm of this convolution by $D_p$.
  (\ref{item:1})
  If $(\lambda_i)\in \ell^p$ then the function
  \begin{equation*}
    f = \sum_i |\lambda_i| 1_{x_iV}
  \end{equation*}
  is in $L^p(G)$ and 
  $\| f\|_{L^p(G)} = |V|\, \|(\lambda_i)\|_{\ell^p}$.
  Convolution with $|F_2|$ is continuous
  so 
  \begin{equation*}
    f*|F_2| 
    = \sum_i |\lambda_i| 1_{x_iV}*|F_2| 
  \end{equation*}
  is in $L^p(G)$.
  Now let us show that the function $1_{x_iV}*|F_2|$ is bigger
  than some constant times $\ell_{x_i}|F_2|$.
  \begin{equation*}
    \int 1_{x_iV}(z) |F_2(z^{-1}y)| dz 
    = \int_V |F_2(z^{-1}x_i^{-1} y)| dz
    \geq C_1 |V| |F_2(x_i^{-1}y)|
  \end{equation*}
  This shows that
  \begin{align*}
    \Big| \sum_i \lambda_i \ell_{x_i} F_2(y)   \Big|
    &\leq     \sum_i |\lambda_i| |F_2(x_i^{-1}y)| \\
    &\leq \frac{1}{C_1|V|}   \sum_i |\lambda_i| 1_{x_iV}*|F_2|(y) \\
    &=  \frac{1}{C_1|V|} f*|F_2|.
  \end{align*}
  Since $f*|F_2|\in L^p(G)$ the sum $\sum_i \lambda_i \ell_{x_i} F_2$ 
  is in $L^p(G)$ with norm
  \begin{align*}
    \Big\| \sum_i \lambda_i \ell_{x_i} F_2(y)   \Big\|_{L^p}
    & \leq \frac{1}{C_1|V|} \|f*|F_2|\, \|_{L^p} \\
    & \leq \frac{D_p}{C_1|V|} \| f\|_{L^p}  \\
    &= \frac{D_p}{C_1}  \| (\lambda_i)\|_{\ell^p}.
  \end{align*}
  This shows the desired continuity.

  (\ref{item:2})
  We need to show that $f(x_i)$ is in $\ell^p$, but this is the same
  as showing that $g = \sum_i |f(x_i)| 1_{x_iV}$ is in $L^p(G)$.
  $f\in L^p(G)*F_2$ so 
  \begin{equation*}
    \sum_i |f(x_i)|1_{x_iV}(y) 
    \leq
    \sum_i |f|*|F_2|(x_i) 1_{x_iV}(y)
    = \int |f(z)| \sum_i |F_2(z^{-1}x_i)| 1_{x_iV}(y) dz
  \end{equation*}
  For each $y$ at most one $i$ adds to this sum, namely the $i$ for which
  $x_i \in yV^{-1}$. Therefore
  \begin{equation*}
      \sum_i |F_2(z^{-1}x_i)| 1_{x_iV}(y)
      \leq \sup_{v\in V} |F_2(z^{-1}yv^{-1})|
      \leq C_2 |F_2(z^{-1}y)|
  \end{equation*}
  by Corollary~\ref{cor:1}. We then get
  \begin{equation*}
    \sum_i |f(x_i)|1_{x_iV}(y) 
    \leq C_2 \int |f(z)| |F_2(z^{-1}y)| dz
    = C_2 |f|*|F_2|(y)
  \end{equation*}
  and finally 
  \begin{equation*}
    \|f(x_i) \|_{\ell^p} \leq \frac{C_2D_p}{|V|} \| f\|_{L^p}.
  \end{equation*}
  (\ref{item:3})
  We have to show that the function
  \begin{equation*}
    \sum_{i} \Big(\int f(x)\psi_i(x)dx \Big) 1_{x_iV} \in L^p(G)
  \end{equation*}
  We get that
  \begin{equation*}
    \Big| 
    \sum_{i} \Big(\int f(x)\psi_i(x)dx \Big) 1_{x_iV}(y)
    \Big|
    \leq  \int |f(x)| \sum_{i\in I} \psi_i(x) 1_{x_iV}(y) dx
  \end{equation*}
  and since 
  \begin{equation*}
    \sum_{i} \psi_i(x)1_{x_iV}(y) 
    \leq \sum_{i} 1_{x_iU}(x)1_{x_iV}(y)
    \leq 1_{U^{-1}V}(x^{-1}y)
  \end{equation*}
  we obtain
  \begin{equation*}
    \Big| 
    \sum_{i} \Big(\int f(x)\psi_i(x)dx \Big) 1_{x_iV}(y)
    \Big|
    \leq
    \int |f(x)|1_{U^{-1}V}(x^{-1}y) dx
    = |f|*1_{U^{-1}V}(y)
  \end{equation*}
  which is in $L^p(G)$ with norm continuously dependent on $f$, i.e.
  \begin{equation*}
    \Big\| \sum_{i} \Big(\int f(x)\psi_i(x)dx \Big) 1_{x_iV}
    \Big\|
    \leq C \| f\|_{L^p}
  \end{equation*}
  for some $C>0$.
\end{proof}

\begin{proposition}
  We can choose a compact neighbourhood $U$, 
  $U$-dense points $\{ x_i \}$ and a partition
  $\psi_i$ of unity with $\supp(\psi_i)\subseteq x_iU$
  such that the operators 
  defined below are invertible with
  continuous inverses
\begin{enumerate}
  \item define $S:L^p(G)*F_2 \to L^p(G)*F_2$
    by
    \begin{equation*}
      Sf = \sum_i f(x_i)\psi_i*F_2
    \end{equation*}
    \label{item:4}
  \item define $T:L^p(G)*F_2 \to L^p(G)*F_2$
    by (with $c_i = \int \psi_i$)
    \begin{equation*}
      Tf = \sum_i c_if(x_i)\ell_{x_i}F_2
    \end{equation*}
    \label{item:5}
  \item define $R:L^p(G)*F_2 \to L^p(G)*F_2$
     by 
    \begin{equation*}
      Rf = \sum_i \Big( \int f(x)\psi_i(x)\, dx\Big)\, \ell_{x_i}F_2
    \end{equation*}
    \label{item:6}
  \end{enumerate}
\end{proposition}

\begin{proof}
  For each neighbourhood of the identity  
  $U$ we can pick $U$-dense points $\{ x_i\}$
  such that $\{ x_i \}$ are $V$-separated for some
  compact neighbourhood of the identity 
  $V$ satisfying $V^2\subseteq U$ (see \cite[Thm 4.2.2]{Rauhut2}).
  Thus we can pick $U$ in order to satisfy the inequality
  in Lemma~\ref{lem:osc} for any $\epsilon$.

  Denote by $D_p$ the $L^p$ operator norm of convolution by
  $F_2$.

  (\ref{item:4})
  Let $f\in L^p(G)*F_2$ and let us look at the difference
  \begin{equation*}
    f(x) - \sum_i f(x_i)\psi_i(x) = \sum_i (f(x)-f(x_i))\psi_i
  \end{equation*}
  For $x\in \supp(\psi_i) \subseteq x_i U$ we get
  \begin{align*}
    |f(x)-f(x_i)|
    &\leq \int |f(z)| |F_2(z^{-1}x) - F_2(z^{-1}x_i)| dz \\
    &= \int |f(z)|  
    \Big| \frac{F_2(z^{-1}x_i)}{ F_2(z^{-1}x)} -1 \Big|
    |F_2(z^{-1}x)| dz
    &\leq \epsilon \int |f(z)| |F_2(z^{-1}x)| dz \\
    &= \epsilon |f|*|F_2|(x)
  \end{align*}
  This means that
  \begin{equation*}
    |f(x) - \sum_i f(x_i)\psi_i(x)|
    \leq 
    \epsilon \sum_i |f|*|F_2|(x) \psi_i(x)
    = \epsilon |f|*|F_2|(x)
  \end{equation*}
  This function is in $L^p(G)$ and so
  \begin{equation*}
    \Big \| f - \sum_i f(x_i)\psi_i \Big\|_{L^p}
    \leq \epsilon D_p \|f \|_{L^p}.
  \end{equation*}
  Convoluting this expression by $F_2$ we get
  \begin{equation*}
    \Big \| f - \sum_i f(x_i)\psi_i*F_2 \Big\|_{L^p}
    \leq \epsilon D_p^2 \|f \|_{L^p}
  \end{equation*}
  So picking $U$ such that $\epsilon < D_p^{-2}$ we obtain 
  an operator $S$ such that 
  $\| I-S\| < 1$ as an operator on $L^p(G)*F_2$.
  Therefore $S$ is invertible.

  (\ref{item:5})
  We will show that $T$ is invertible using its difference from
  the operator $S$.
  \begin{align*}
    |Tf(x) - Sf(x)|
    &= \Big| \sum_i f(x_i) (\psi_i*F_2(x) - c_i F_2(x_i^{-1}x)) \Big| \\
    &\leq \sum_i |f(x_i)| |\psi_i*F_2(x) - c_i F_2(x_i^{-1}x)|
  \end{align*}
  Look at 
  \begin{align*}
    |\psi_i*F_2(x) - c_i F_2(x_i^{-1}x)|
    &= \Big| 
    \int \psi_i(z) (F_2(z^{-1}x) -F_2(x_i^{-1}x) dz
    \Big| \\
    &\leq \int \psi_i(z) |F_2(z^{-1}x) -F_2(x_i^{-1}x)| dz \\
    &\leq \epsilon \int \psi_i(z) |F_2(z^{-1}x)| dz \\
    &= \epsilon \psi_i*|F_2|(x)
  \end{align*}
  Then we have 
  \begin{equation*}
    |Tf(x) - Sf(x)|
    \leq \epsilon \sum_i |f(x_i)| \psi_i*|F_2|(x)
  \end{equation*}
  This is a function in $L^p(G)$ and the norm is
  \begin{equation*}
    \|Tf-Sf \|_{L^p} 
    \leq \epsilon D_p \Big\| \sum_i |f(x_i)| \psi_i \Big\|_{L^p}
    \leq \epsilon^2 D_p^2  \| f\|_{L^p}
  \end{equation*}
  This means that
  $\| I-T\| \leq \| I-S\| + \| S-T\| \leq \epsilon D_p(1 + \epsilon D_p)$
  and if we pick $U$ such that this norm is less than $1$ we get
  that $T$ is invertible.

  (\ref{item:6})
  Can be proven by estimates as above.

\end{proof}

This now means that any $f\in L^p(G)*F_2$ can be reconstructed
in the following ways (we only write this up for the operator $T$)
\begin{align*}
  f &= \sum_i c_i f(x_i) T^{-1}(\ell_{x_i}F_2)
  \intertext{or}
  f &= \sum_i c_i (T^{-1}f)(x_i)\ell_{x_i}F_2
\end{align*}

The first representation in turn means that a 
$v'\in \Co_{\mathcal{H}_n^\infty}^u L^p(G)$
can be reconstructed from the samples of its
voice transform, i.e.
\begin{equation*}
  v' = \sum_i c_i V_u(v')(x_i) V_u^{-1}T^{-1}\ell_{x_i}V_u(u)
\end{equation*}

\begin{remark}
  In \cite{FG1} it is concluded that the elements
  $v_i = V_u^{-1}T^{-1}\ell_{x_i}V_u(u)$ are in the space 
  $\mathcal{H}_m^1$. We cannot obtain such a result, since
  this space does not exist in general. We however claim, 
  that the $v_i$ are in all the
  coorbit spaces $\Co_{\mathcal{H}_n^\infty}^u L^p(G)$
  for $p>1$ and so the same vectors can be used in all 
  situations.
\end{remark}

\subsection{Bandlimited functions}

\label{sec:bandl-funct}
\noindent
In this section we give a examples of coorbit spaces corresponding to
cyclic representations that are far from irreducible.  The
examples are almost tautological, the representation space is the
coorbit space itself. 
This makes it possible to treat sampling of bandlimited
functions using the coorbit theory. 

Let $G$ be a locally compact group
and let $K\subset G$ be a compact subgroup.  We
assume that $(G,K)$ is a Gelfand pair, i.e. that the algebra
$L^1(G /K)^K$ of left-invariant functions with convolution
is commutative. Thus
\begin{equation}\label{eq-DirectInt}
  (\ell,L^2(G/K))
  \simeq 
  \int^\oplus_{\Lambda }(\pi_\lambda ,\mathcal{H}_\lambda )\, d\mu (\lambda )
\end{equation}
with multiplicity one. Here $\ell$ stands for the left regular
representation $[\ell (x)f] (y)=f(x^{-1}y)$ and $\Lambda$ is a
measurable subset of $\widehat{G}$. 
Those spaces are sometimes called
\textit{commutative space} \cite{joe07}.  For details in the following
arguments we refer to \cite{joe07} or \cite{os07} for the case of
Riemannian symmetric spaces of the non-compact type.  Because of the
multiplicity one assumption it follows that $\dim \mathcal{H}_\lambda^K=1$ for
almost all $\lambda\in\Lambda$. Assume that we can choose a measurable
field $\Lambda\ni\lambda \mapsto s_\lambda\in \mathcal{H}_\lambda^K$ such that
$\|s_\lambda \|=1$ for almost all $\lambda$. Note that $s_\lambda$ is then
determined up to a constant $c_\lambda$, $|c_\lambda |=1$. The
function $x\mapsto \pi_\lambda (x)s_\lambda$ is right $K$-invariant
and can therefore be viewed as a function on $G/K$. From now on
we will indentify functions on $G/K$ with right $K$-invariant
functions on $G$.  For $f\in L^1(G/K)$ let
$$\pi_\lambda (f)s_\lambda =\int_G f(x)\pi_\lambda (x)s_\lambda\, dx =
\int_{G/K} f(x)\pi_\lambda (x)s_\lambda \, dx \, .$$ The map
$\mathcal{F} : L^2(G/K)\to \int_\Lambda^\oplus \mathcal{H}_\lambda\, d\mu
(\lambda )$ defined by
\begin{equation}\label{eq-FourierTr}
  \mathcal{F} (f)(\lambda )=\widehat{f}(\lambda ):=\pi_\lambda (f)s_\lambda
\end{equation}
the \textit{vector valued Fourier transform} on $G/K$, or simply, as we
will not be considering another Fourier transform here, the Fourier
transform on $G/K$. First $\mathcal{F}(f)$ is only defined for $f\in
L^1(G/K)\cap L^2(G/K)$ but the following argument shows that
$\|f\|=\|\mathcal{F}(f)\|$ and hence $\mathcal{F}$ extends to a
unitary isomorphism giving a ``concrete'' realization of the
isomorphism (\ref{eq-DirectInt}).

As $\pi_\lambda (\ell (x)f)= \pi_\lambda (x)\pi_\lambda (f)$ and
$\pi_\lambda (f*g)=\pi_\lambda (f)\pi_\lambda (g)$ ($x\in G$, $f,g\in
L^1(G/K)\cap L^2(G/K)$) it follows that:
\begin{lemma}\label{le-1.1} Let $x\in G$ and $f,g\in L^1 (G/K)\cap
L^2(G/K)$. Then 
$\mathcal{F} (\ell (x)f)(\lambda)=\pi_\lambda (x)\widehat{f}(\lambda)$
  and $\widehat{f*g}(\lambda )=\pi_\lambda (f)\widehat{g}(\lambda )$.
\end{lemma}

We also have
\begin{eqnarray*}
  \mathrm{Tr} (\pi_\lambda (f))&=&(\widehat{f}( \lambda ),s_\lambda )_\lambda \\
  &=&\int_{G/K}f(x)(\pi_\lambda (x)s_\lambda ,s_\lambda )_\lambda \, dx\\
  &=&\int_{G/K} f(x)\overline{\varphi_\lambda (x)}\, dx
\end{eqnarray*}
where $\varphi_\lambda$ is the \textit{spherical function}
$\varphi_\lambda (x)=(s_\lambda, \pi_\lambda (x)s_\lambda )_\lambda $.
The spherical function is independent of the choice of $s_\lambda$.
Using the well known theory of the spherical Fourier transform
\cite{joe07} or the abstract Plancherel formula it follows that for 
$f\in C_c(G/K)$ we can choose the measure $\mu$ on $\Lambda$ such that
$$f(eK)=\int_\lambda (\widehat{f}(\lambda ),s_\lambda )_\lambda \, d\mu (\lambda )\, .$$
Define $f^*(x)=\Delta (x^{-1})\overline{f(x^{-1})}$. Then $\pi_\lambda
(f^*)=\pi_\lambda (f)^*$. Thus $\widehat{f^* * f}= \pi_\lambda
(f)^*\widehat{f}(\lambda)$ and hence
\begin{equation*}
    \|f\|_2^2 = \int_\Lambda \|\widehat{f} (\lambda ) \|^2_\lambda \, d\mu (\lambda )\, .
\end{equation*}

\begin{lemma}\label{le-inversion} Let $f\in L^1(G/K)\cap L^2(G/K)$ be such that
  $\lambda \mapsto \|\widehat{f}(\lambda )\| $ is integrable. Then $f$
  is continuous and
$$f(x)= \int_\Lambda f*\varphi_\lambda (x)\, d\mu (\lambda )=
\int_\Lambda (\widehat{f} (\lambda ),\pi_\lambda (x )s_\lambda
)_\lambda \, d\mu (\lambda ) \, .$$
\end{lemma}
\begin{proof} Note first that
  \begin{equation}\label{eq-norm}
    | (\widehat{f}(\lambda ),\pi_\lambda (x) s_\lambda)_\lambda |
    \le \|\widehat{f}(\lambda )\|_\lambda \|\pi_\lambda (x) s_\lambda\|_\lambda
    = \|\widehat{f}(\lambda )\|_\lambda \, .
  \end{equation}
  Using (\ref{eq-FourierTr}) we get
  \begin{eqnarray*}
    f(x )&=&\ell(x^{-1})f(e)\\
    &=&\int_\Lambda (\mathcal{F}(\ell (x^{-1})f)(\lambda )
    ,s_\lambda )_\lambda \, d\mu (\lambda )\\
    &=&\int_\Lambda (\pi_\lambda (x^{-1})\widehat{f} (\lambda ),s_\lambda )_\lambda \, d\mu (\lambda )\\
    &=&\int_\Lambda (\widehat{f} (\lambda ),\pi_\lambda (x )s_\lambda )_\lambda \, d\mu (\lambda )
    \\
    &=&\int_\Lambda \int_{G/K} f(y)(s_\lambda ,\pi_\lambda (y^{-1}x)s_\lambda )_\lambda\, dy d\mu (\lambda )\\
    &=&\int_\Lambda f*\varphi_\lambda (x)\, d\mu (\lambda )\, .
  \end{eqnarray*}
  The rest follows now from (\ref{eq-norm}).
\end{proof}
\begin{definition} Let $\Omega\subset \Lambda$ be measurable with
  finite measure.  The function $f\in L^2(G/K)$ is called
  $\Omega$-bandlimited if $\supp (\widehat{f}) \subseteq \Omega$. We
  denote the space of $\Omega$-bandlimited functions by $L^2_\Omega
  (G/K)$.
\end{definition}

From now on we will always assume that $\Omega\subset\Lambda$ has
finite measure. Define a section $\widehat{\psi}\in
\int^\oplus_\Lambda (\pi_\lambda, \mathcal{H}_\lambda )\, d\mu (\lambda )$ by
\begin{equation}\label{eq-psi}
  \widehat{\psi} (\lambda )=\chi_{\Omega }(\lambda )s_\lambda \, .
\end{equation}
where $\chi_\Omega$ is the characteristic function on $\Omega$.
Note that $\|\widehat{\psi}\|=\sqrt{\mu (\Omega )}<\infty$.  Let $\psi
:=\mathcal{F}^{-1}(\widehat{\psi })$. Then $\psi\in L^2_\Omega (G/K)$ and
$\|\psi\|_2=\sqrt{\mu (\Omega )}$. By Lemma \ref{le-inversion} and the
fact that $\overline{\varphi_\lambda (x^{-1})}=\varphi_\lambda (x)$ we
get:

\begin{lemma} We have $\psi (x)=\int_\Omega \varphi_\lambda (x)\, d\mu
  (\lambda )$ and $\psi (x)=\overline{\psi (x^{-1})}$.
\end{lemma}

\begin{lemma}\label{le-invariant} Let the notation be as above. Then
  the following holds:
  \begin{enumerate}
  \item $L_\Omega^2(G/K)$ is a $G$-invariant Hilbert space. \label{item:7}
  \item The orthogonal projection onto $L^2_\Omega (G/K)$ is given by
    $f\mapsto f*\psi$. In particular $\psi *\psi =\psi$ and $\psi$ is
    cyclic in $L^2_\Omega (G/K)$.\label{item:8}
  \item $L_\Omega^2(G/K)\subset C(G/K)$. \label{item:9}
  \item If $f\in L^2_\Omega (G/K)$ then $|f(x)|\le \sqrt{\mu (\Omega
      )}\|f\|_2$.  In particular, the point evaluation map
    $\mathrm{ev}_x : L^2(G/K) \to \mathbb{C}$, $f\mapsto f(x)$, $x\in
    G/K$ is continuous. \label{item:10}
  \end{enumerate}
\end{lemma}

\begin{proof} (\ref{item:7}) and (\ref{item:8}) 
  follow from Lemma~\ref{le-1.1} $\mathcal{F}
  (\ell(x)f)= \pi_\lambda (x)\widehat{f}(\lambda )$ and
$$\mathcal{F} (f*\psi )(\lambda )=\pi_\lambda (f)\mathcal{F} (\psi )(\lambda )
=\chi_\Omega (\lambda )\pi_\lambda (f)s_\lambda = \chi_\Omega (\lambda
) \mathcal{F} (f) (\lambda )\, .$$ For (\ref{item:9}) 
we note that if $f\in L^2_\Omega
(G/K)$ then $f=f*\psi$ and hence $f$ is continuous. It follows that
the evaluation map $\mathrm{ev}_x$, $x\in G/K$ is well defined.
Furthermore, $|f(x)|=|f*\psi (x)|\le \|f\|_2\|\psi \|_2$ and (\ref{item:10})
follows.
\end{proof}

\begin{remark} This shows that the spaces of bandlimited functions are
  natural spaces to study sampling. In particular each of those spaces
  is a reproducing kernel Hilbert space with reproducing kernel
  $K(x,y)=\overline{\psi (y^{-1}x)}$. 
\end{remark}

\begin{theorem} Let $\Omega\subset \Lambda$ be measurable with finite
  measure.  Let $S=L^2_\Omega (G/K)$.  Then $L^2_\Omega
  (G/K)=\mathrm{Co}_{S}^\psi (L^2(G/K))$ is a Coorbit space.
\end{theorem}

\begin{proof} We have for $f\in L^2(G/K)$
  $$
  V_\psi (f)(x)= (f,\lambda (x)\psi )=f*\psi (x)\, .$$ The claim follows
  therefore from Lemma \ref{le-invariant}.
\end{proof}

\begin{remark} If $G$ is a Lie group, then we can replace $S=L^2_\Omega
  (G/K)$ by the space $S=L^2_\Omega (G/K)^\infty$ of smooth vectors,
  i.e., all $f\in L^2_\Omega (G/K)$ such that $Df\in L^2_\Omega (G/K)$
  for all invariant differential operators $D$.
\end{remark}

\begin{remark}
  If we can choose a topology on $\Lambda$ such that
  $\lambda\mapsto \pi_\lambda(x)s_\lambda$ is continuous and $\Omega$
  is chosen to be compact, then
  the function $\psi$ cannot be integrable. If it were this would mean
  that $\widehat{\psi}$ is continuous, which is not the case.
  Therefore we have another example for which integrability is not
  necessary.  
\end{remark}

\begin{example}[$G$ abelian] If $G$ is abelian then
$$L^2(G)\simeq L^2(\widehat{G},d\mu )\simeq \int^\oplus_{\widehat{G}}(\chi ,\mathbb{C})\, d\mu (\chi)$$
where $\widehat{G}$ is the dual group consisting of all continuous characters
$\chi : G\to \mathbb{T}$ and $d\mu $ is the dual Haar measure on
$\widehat{G}$. In this case we can choose $s_\chi =1$ for all
$\chi\in\widehat{G}$. The Fourier transform is given by
$$\widehat{f}(\chi )=\int_G f(x)\chi (x)\, dx\, .$$
Note, that the usual conjugation is missing. The function $\psi$ is
$$\psi (x)=\int_\Omega \overline{\chi (x)}\, d\mu (\chi )=
\int_\Omega \chi (x^{-1} ) \, d\mu (\chi )\, .$$

If $G=\mathbb{R}$ and $\Omega=[-R/2,R/2]$ then
$$\psi (x)=\int_{-R/2}^{R/2} e^{-2\pi i\lambda x}\, d\lambda
=\frac{\sin (R \pi x)}{\pi x}=R \, \mathrm{sinc}(R\pi x)$$ as
usually. The Shannon sampling theorem is then
$$f(x)=\sum_{n=-\infty}^\infty f(n/R)\frac{\sin (R\pi ( x - n/R )}{\pi (x - n/R)}\, .$$

In general, if there exists a discrete subgroup
$\widehat{\Gamma}\subset \widehat{G}$ such that $\Omega \simeq
\widehat{G}/\widehat{\Gamma}$, ie, $\Omega +\widehat{G}$ is a
measurable tiling of $\widehat{G}$ then we have with the same proof as
in the classical case:
$$f(y)=\sum_{x\in \widehat{\Gamma}^\perp }f (x)\psi (x^{-1}y)$$
uniformly. Here
$$\widehat{\Gamma}^\perp =\{x\in G\mid (\forall \gamma\in \widehat{\Gamma})\,\, \gamma (x)=1\}\, .$$
The interesting case is therefore the irregular sampling. That
discussion has been started in \cite{Fei92}.
\end{example}
\begin{example}[Compact groups] Let $K$ be a compact group. Set
  $G=K\times K$ and identify $K$ with the diagonal subgroup $K\simeq
  \{(k,k)\in G\mid k\in K\}$. The space $G/K$ can be identified with
  $K$ by $(a,b)K\mapsto ab^{-1}$. The left regular action on $G/K$ is
  then transformed to the left-right action of $K$ on $L^2(K)$, i.e.,
  $\ell_K\otimes r_K(a,b)f(x) =f(a^{-1}xb)$.
  For an irreducible representation $(\pi, V_\pi )$ of $K$ define a
  representation of $G$ on $\mathrm{End}(V_\pi )$ by
$$\widetilde{\pi }(a,b)T=\pi (a)\circ T\circ \pi (b^{-1})\, .$$
Then
$$L^2(K)=L^2(G/K)\simeq \bigoplus_{\pi\in \widehat{G}}(\widetilde{\pi },\mathrm{End}(V_\pi ))\,.$$
We choose $s_\pi=\mathrm{I}_{V_\pi}$. Then $\widehat{f} (\pi ) =\pi
(f^\vee )I_{V_\pi}$ where $f^\vee (x)=f(x^{-1})$. The inversion formula
reads
$$f(x)=\sum_{\pi }d (\pi )\mathrm{Tr} (\pi (x)\widehat{f}(\pi ))\, .$$
A subset $\Omega\subset \widehat{G}$ is of finite measure if and only
if $\Omega$ is finite. In that case we have
$$\psi (x)=\sum_{\pi\in \Omega} d(\pi )\mathrm{Tr} (\pi (x)I_{V_\pi})$$
and $\widehat{\psi}(\pi )=I_{V_\pi}$. Furthermore
$$L^2_\Omega (G/K)=\Big\{\sum_{\pi \in \Omega}d(\pi )\mathrm{Tr} (\pi (x)T_\pi )\Big| (\forall \pi \in \Omega )\, \, T_\pi
\in \mathrm{End}(V_\pi)\Big\}\, .$$
\end{example}
\begin{example}[Riemannian symmetric spaces of the noncompact type]
  The last example we would like to mention is the case of Riemannian
  symmetric spaces $G/K$ of the non-compact type, which always comes
  from a Gelfand pair $(G,K)$. We will not describe this
  case here in detail as it would require us to introduce several pages of new
  notation. We would however like to point out the article by Pesenson
  \cite{pes07} where irregular sampling is discussed. Note, that
  according to \cite{kst}, see also \cite{os07} every function $f$ in
  $L^2_\Omega (G/K)$ extends to a holomorphic function on a complex
  neighborhood of $G/K$ and hence $f$ is determined on any infinite
  set with a limit point.
\end{example}

\bibliographystyle{plain}
\bibliography{coorbit}
\end{document}